\documentclass{amsart}
\usepackage{amssymb}
\usepackage{amsmath}

\newtheorem{theorem}{Theorem}[section]
\newtheorem{lemma}[theorem]{Lemma}
\newtheorem{proposition}[theorem]{Proposition}
\theoremstyle{remark}
\newtheorem*{remark}{Remark}
\newtheorem*{remarks}{Remarks}
\numberwithin{equation}{section}
\begin{document}

\title{Generators of relations for annihilating fields}
\author{Mirko Primc}
\address{Dept. of Math., Univ. of Zagreb, Bijeni\v{c}ka 30, 10000 Zagreb, Croatia}
\email{primc@math.hr}
\thanks{Partially supported  by the Ministry of Science and Technology of the
Republic of Croatia, grant 037002.} 

\begin{abstract}
For an untwisted affine Kac-Moody Lie algebra $\tilde{\mathfrak
g}$, and a given positive integer level $k$, vertex operators
$x(z)=\sum x(n)z^{-n-1}$, $x\in\mathfrak g$, generate a vertex
operator algebra $V$. For the maximal root $\theta$ and a root
vector $x_\theta$ of the corresponding finite-dimensional
$\mathfrak g$, the field $x_\theta(z)^{k+1}$ generates all
annihilating fields of level $k$ standard $\tilde{\mathfrak
g}$-modules. In this paper we study the kernel of the normal order
product map $r(z)\otimes Y(v,z)\mapsto :r(z) Y(v,z):$ for $v\in V$
and $r(z)$ in the space of annihilating fields generated by the
action of $\tfrac{d}{dz}$ and $\mathfrak g$ on
$x_\theta(z)^{k+1}$. We call the elements of this kernel the
relations for annihilating fields, and the main result is that
this kernel is generated, in certain sense, by the relation
$x_\theta(z)\tfrac{d}{dz}(x_\theta(z)^{k+1})=
(k+1)x_\theta(z)^{k+1}\tfrac{d}{dz}x_\theta(z)$. This study is
motivated by Lepowsky-Wilson's approach to combinatorial
Rogers-Ramanujan type identities, and many ideas used here stem
from a joint work with Arne Meurman.
\end{abstract}
\maketitle

\section{Introduction}

J.~Lepowsky and R.~L.~Wilson (cf. \cite{LW}) initiated the
approach to combinatorial Rogers-Ramanujan type identities via the
vertex operator constructions of representations of affine Lie
algebras. By using Lepowsky-Wilson's approach, in \cite{MP1} and
\cite{FKLMM} new series of combinatorial Rogers-Ramanujan type
identities are obtained by constructing certain combinatorial
bases of standard $\mathfrak{sl}(2,\mathbb C)\widetilde{}\
$-modules. In both papers the starting point is a PBW spanning set
of a standard module $L(\Lambda)$ of level $k$, which is then
reduced to a basis by using the relation
\begin{equation}\label{E:1.1}
x_\theta(z)^{k+1}=0\quad\text{on}\quad L(\Lambda).
\end{equation}

The relation (\ref{E:1.1}) holds in general for any level $k$
standard module of any untwisted affine Kac-Moody Lie algebra
$\tilde{\mathfrak g}$. Here $\theta$ is the maximal root of the
corresponding finite-dimensional Lie algebra $\mathfrak g$,
$x_\theta$ is a corresponding root vector and
$x_\theta(z)=\sum_{n\in\mathbb Z}x_\theta(n)z^{-n-1}$ is a formal
Laurent series with coefficients $x_\theta(n)$ in the affine Lie
algebra $\tilde{\mathfrak g}$ (see Section~2 for notation). The
relation (\ref{E:1.1}) was studied before in \cite{LP}, but it is
best understood if $x_\theta(z)^{k+1}$ is viewed as the vertex
operator $Y(x_\theta(-1)^{k+1}\mathbf 1,z)$ associated to a
singular vector $x_\theta(-1)^{k+1}\mathbf 1$ in a vertex operator
algebra $V=N(k\Lambda_0)$, i.e., as an annihilating field of
standard modules of level $k$ (cf. \cite{FLM}, \cite{DL},
\cite{Li1}, \cite{MP1}).

After a PBW spanning set is reduced to a basis, it remains to
prove its linear independence. In a joint work with Arne Meurman
\cite{MP1} a sort of Gr\"obner base theory is used for proving
linear independence. The main ingredient in the proof is the use
of relations among annihilating fields of level $k$ standard
modules
$$
Y(r,z)=\sum_{n\in\mathbb Z}r_nz^{-n-1},\qquad r\in \bar R {\mathbf
1}=\mathbb C[L_{-1}]U(\mathfrak g)x_\theta(-1)^{k+1}\mathbf 1,
$$
``generated'' by the obvious relation
\begin{equation}\label{E:1.2}
x_\theta(z)\tfrac{d}{dz}(x_\theta(z)^{k+1})=
(k+1)x_\theta(z)^{k+1}\tfrac{d}{dz}x_\theta(z).
\end{equation}
The relation (\ref{E:1.2}) arises by writing the vector
$x_{\theta}(-2)x_\theta(-1)^{k+1}\mathbf 1$ in two different ways,
thus obtaining two different expressions for the same field, and
the representation theory is used to generate ``enough'' of other
relations for annihilating fields $Y(r,z)$, $r\in \bar R {\mathbf
1}$. In \cite{MP2} a similar technique is used to construct a
combinatorial basis of the basic $\mathfrak{sl}(3,\mathbb
C)\widetilde{}\ $-module.

By following ideas developed in \cite{MP1} and \cite{MP2}, in
\cite{P2} a general construction of such relations for
annihilating fields is given by using vertex operators $Y(q,z)$
associated with elements $q$ in the kernel of the map
$$
\Phi \colon \bar R {\mathbf 1}\otimes N(k\Lambda_0)\rightarrow
N(k\Lambda_0),\quad \Phi (u\otimes v)=u_{-1}v.
$$
Due to the state-field correspondence, we shall call the elements
in $\ker\Phi$ {\it relations for annihilating fields}. By using
these relations the problem of constructing combinatorial bases of
standard modules is split into a ``combinatorial part of the
problem'' (consisting of counting numbers $N(n)$ of certain
colored partitions of $n$) and a ``representation theory part of
the problem '' (consisting of constructing certain subspaces $Q$
in $\ker\Phi$ such that $\dim Q(n)=N(n)$, where $Q(n)$ is a space
of coefficients of vertex operators $Y(q,z)=\sum
q(n)z^{-n-\text{wt\,}q}\,$ associated with vectors in $Q\,$). In
this paper only a part of the ``representation theory part of the
problem '' is studied, the main objective being to understand the
structure of $\ker\Phi$ for general untwisted affine Kac-Moody Lie
algebra $\tilde{\mathfrak g}$.

The starting point in this paper are Theorems~\ref{T:3.1} and
\ref{T:3.2}, showing that $\ker\Phi$ has certain $(\mathbf C
L_{-1}\ltimes\tilde{\mathfrak g}_{\leq 0})$-module structure and
that
\begin{equation}\label{E:1.3}
\ker\Phi\cong\ker\Psi
\end{equation}
as $(\mathbf C L_{-1}\ltimes\tilde{\mathfrak g}_{\leq
0})$-modules, where
$$
\Psi \colon U(\tilde{\mathfrak g})\otimes_{U(\tilde{\mathfrak
g}_{\geq 0}+\mathbb C c)} \bar R {\mathbf 1} \to
N(k\Lambda_0),\qquad u\otimes w \mapsto uw.
$$
This makes it possible to use the theory of induced
representations of $\tilde{\mathfrak g}$, in particular
Garland-Lepowsky's resolution of a standard module in terms of
generalized Verma modules. The main results are
Theorem~\ref{T:5.4}, saying that the $(\mathbf C
L_{-1}\ltimes\tilde{\mathfrak g}_{\leq 0})$-module $\ker\Phi$ is
generated by certain (singular) vectors, and for combinatorial
arguments more convenient Theorem~\ref{T:6.2}, saying that the
$(\mathbf C L_{-1}\ltimes\tilde{\mathfrak g})$-module $\ker\Psi$
is generated by the vector
$$
x_{\theta}(-2)\otimes x_\theta(-1)^{k+1}\mathbf 1\mathbf
-x_{\theta}(-1)\otimes x_{\theta}(-2)x_\theta(-1)^{k}\mathbf 1.
$$
Via (\ref{E:1.3}) this vector can be seen in the vertex operator
algebra $N(k\Lambda_0)\otimes N(k\Lambda_0)$, and the vertex
operator associated to this vector is
\begin{equation}\label{E:1.4}
x_\theta(z)^{k+1}\otimes\tfrac{d}{dz}x_\theta(z)-
\tfrac{1}{k+1}\,\tfrac{d}{dz}(x_\theta(z)^{k+1})\otimes
x_\theta(z).
\end{equation}
The obvious relation (\ref{E:1.2}) arises by taking from
(\ref{E:1.4}) normal order products. So, in a sense, the relation
(\ref{E:1.2}) generates all relations for annihilating fields of
standard $\tilde{\mathfrak g}$-modules of level $k$.

Some ideas worked out in this paper come from collaboration with
Arne Meurman for many years, and I thank Arne Meurman for his
implicit contribution to this work. I also thank Jim Lepowsky and
Ivica Siladi\' c for numerous stimulating discussions.

\section{Vertex algebras for affine Lie algebras}

Let ${\mathfrak g}$ be a simple complex Lie algebra, $\mathfrak h$
a Cartan subalgebra of ${\mathfrak g}$ and $\langle \ , \ \rangle$
a symmetric invariant bilinear form on ${\mathfrak g}$. Via this
form we identify $\mathfrak h$ and $\mathfrak h^*$ and we assume
that $\langle \theta , \theta \rangle=2$ for the maximal root
$\theta$ (with respect to some fixed basis of the root system).
Set
$$
\hat{\mathfrak g} =\coprod_{j\in\mathbb Z}{\mathfrak g}\otimes
t^{j}+\mathbb C c, \qquad \tilde{\mathfrak g}=\hat{\mathfrak
g}+\mathbb C d.
$$
Then $\tilde{\mathfrak g}$ is the associated untwisted affine Lie
algebra (cf. \cite{K}) with the commutator
$$
[x(i),y(j)]=[x,y](i+j)+i\delta_{i+j,0}\langle x,y\rangle c.
$$
Here, as usual, $x(i)=x\otimes t^{i}$ for $x\in{\mathfrak g}$ and
$i\in\mathbb Z$, $c$ is the canonical central element, and
$[d,x(i)]=ix(i)$. Sometimes we shall denote ${\mathfrak g}\otimes
t^{j}$ by ${\mathfrak g}( {j})$. We identify ${\mathfrak g}$ and
${\mathfrak g}(0)$. Set
$$
\tilde{\mathfrak g}_{<0} =\coprod_{j<0}{\mathfrak g}\otimes
t^{j},\qquad \tilde{\mathfrak g}_{\leq 0} =\coprod_{j\leq
0}{\mathfrak g}\otimes t^{j}+\mathbb C d,\qquad\tilde{\mathfrak
g}_{\geq 0} =\coprod_{j\geq 0}{\mathfrak g}\otimes t^{j}+\mathbb C
d.
$$

For $k\in\mathbb C$ denote by $\mathbb C v_k$ the one-dimensional
$(\tilde{\mathfrak g}_{\geq 0}+\mathbb C c)$-module on which
$\tilde{\mathfrak g}_{\geq 0}$ acts trivially and $c$ as the
multiplication by $k$. The untwisted affine Lie algebra
$\tilde{\mathfrak g}$ gives rise to the vertex operator algebra
(see \cite{FZ} and \cite{Li1}, we use the notation from
\cite{MP1})
$$
N(k\Lambda_0)=U(\tilde{\mathfrak
g})\otimes_{U(\tilde{\mathfrak g}_{\geq 0}+\mathbb C c)}\mathbb C
v_k
$$
for level $k\neq -g^\vee$, where $g^\vee$ is the dual Coxeter
number of ${\mathfrak g}$; it is generated by the fields
$$
x(z)=\sum_{n\in\mathbb Z}x_nz^{-n-1},\qquad x\in{\mathfrak g},
$$
where we set $x_n=x(n)$ for $x\in{\mathfrak g}$. As usual, we
shall write $Y(v,z)=\sum_{n\in\mathbb Z}v_nz^{-n-1}$ for the
vertex operator (field) associated with a vector $v\in
N(k\Lambda_0)$. From now on we shall fix the level $k\neq
-g^\vee$, and we shall sometimes denote by $V$ the vertex operator
algebra $N(k\Lambda_0)$.

Recall that $N(k\Lambda_0)\otimes N(k\Lambda_0)$ is a vertex
operator algebra; fields are defined by $Y(a\otimes
b,z)=Y(a,z)\otimes Y(b,z)$, the conformal vector is $\omega\otimes
\mathbf 1+\mathbf 1\otimes \omega$ (cf.  \cite{FHL}). In
particular, the derivation $D=L_{-1}$ is given by $D\otimes
1+1\otimes D$, the degree operator $-d=L_0$ is given by
$L_0\otimes 1+1\otimes L_0$, and we have the action of ${\mathfrak
g}={\mathfrak g}(0)$ given by $x\otimes 1+1\otimes x$. Set
\begin{equation*}
\Phi \colon N(k\Lambda_0)\otimes N(k\Lambda_0)\rightarrow
N(k\Lambda_0),\qquad \Phi(a\otimes b)=a_{-1}b.
\end{equation*}
It is easy to see that $\Phi$ intertwines the actions of $L_{-1}$,
$L_{0}$ and ${\mathfrak g}(0)$, so that $\ker \Phi$ is also
invariant for the actions of these operators. If $\sigma$ is an
automorphism of $V$, then clearly $\Phi$ intertwines the actions
of $\sigma$ on $V\otimes V$ and $V$.

\section{Tensor products and induced representations}

The vertex operator algebra $V=N(k\Lambda_0)$ is an induced
$\tilde{\mathfrak g}$-module, so by restriction it is a
$\tilde{\mathfrak g}_{<0}$-module, where $x_j=x(j)$,
$x\in\mathfrak g$, $j<0$, acts on
\begin{equation}\label{E:3.1}
V=U(\tilde{\mathfrak g})\otimes_{U(\tilde{\mathfrak g}_{\geq
0}+\mathbb C c)}\mathbb C v_k\cong U(\tilde{\mathfrak g}_{<0})
\end{equation}
by the left multiplication with $x_j$. Since in general
\begin{equation}\label{E:3.2}
(Du)_n=-nu_{n-1}
\end{equation}
for $u\in V$ and $n\in\mathbb Z$, we have
\begin{equation}\label{E:3.3}
x_{-1-j}=(D^{(j\,)}x(-1)\mathbf 1)_{-1}=(x_{-1-j}\mathbf 1)_{-1}
\end{equation}
for $x\in\mathfrak g$ and $j\in\mathbb Z_{\geq 0}$, with notation
$D^{(j\,)}=D^{j}/j!\,$.

The vertex operator algebra $V$ has a Lie algebra structure with
the commutator
\begin{equation}\label{E:3.4}
[u,v]=u_{-1}v-v_{-1}u=\sum_{n\geq 0}(-1)^nD^{(n+1)}(u_nv).
\end{equation}
This Lie algebra structure is transported from the Lie algebra
structure on the Lie algebra $\mathcal L_-(V)$ associated with the
vertex Lie algebra $V$ (cf. \cite{FF}, \cite{Li2}), via the vector
space isomorphism
\begin{equation}\label{E:3.5}
V\to\mathcal L_-(V),\qquad u\mapsto u_{-1}
\end{equation}
(cf. \cite{DLM}, or (4.14) in \cite{P1} with $U=V$). By
(\ref{E:3.3}) the restriction of (\ref{E:3.5}) gives a vector
space isomorphism
\begin{equation}\label{E:3.6}
\tilde{\mathfrak g}_{<0}\mathbf 1\to\tilde{\mathfrak
g}_{<0},\qquad u\mapsto u_{-1},
\end{equation}
and the adjoint action of $\tilde{\mathfrak g}_{<0}\mathbf 1$ on
$V$ with the commutator (\ref{E:3.4}) is transported from the
adjoint action of the subalgebra $\tilde{\mathfrak g}_{<0}$ on the
Lie algebra $\mathcal L_-(V)$. So we have the ``adjoint'' action
$u_{-1}\colon v\mapsto [u,v]$ of the Lie algebra $\tilde{\mathfrak
g}_{<0}$ on the vertex operator algebra $V$, we shall denote it by
$V_{\text{ad}}$.

Since $L_{-1}$, $L_{0}$ and $y_0$, $y\in {\mathfrak g}$, are
derivations of the product $u_{-1}v$ in $V$, they are also
derivations of the bracket $[u,v]$, and we can extend the adjoint
action of the Lie algebra $V$ to the action of the Lie algebra
$\big(\mathbf C L_{-1}+\mathbf C L_{0}+{\mathfrak
g}(0)\big)\ltimes V$. Since $\tilde{\mathfrak g}_{<0}\mathbf 1$ is
invariant for $L_{-1}$, $L_{0}$ and $y_0$, $y\in {\mathfrak g}$,
and
$$
\big(\mathbf C L_{-1}+\mathbf C L_{0}+{\mathfrak g}(0)\big)\ltimes
\tilde{\mathfrak g}_{<0}\mathbf 1\cong \mathbf C
L_{-1}\ltimes\tilde{\mathfrak g}_{\leq 0},
$$
we can extend the ``adjoint'' action of the Lie algebra
$\tilde{\mathfrak g}_{<0}$ on the vertex operator algebra $V$ to
the ``adjoint'' action of the Lie algebra $\mathbf C
L_{-1}\ltimes\tilde{\mathfrak g}_{\leq 0}$ on $V$, we shall also
denote it by $V_{\text{ad}}$.

Let $W\subset V$ be a $\tilde{\mathfrak g}_{\geq 0}$-submodule
invariant for the action of $D=L_{-1}$. Then the right hand side
of (\ref{E:3.4}), and (\ref{E:3.2}), imply that $W$ is invariant
for the ``adjoint'' action of $\tilde{\mathfrak g}_{<0}$. By
assumption $W$ is invariant for  $L_{-1}$, $L_{0}$ and $\mathfrak
g(0)$, so on $W$ we have the ``adjoint'' action of Lie algebra
$\mathbf C L_{-1}\ltimes\tilde{\mathfrak g}_{\leq 0}$, we shall
denote it by $W_{\text{ad}}$.

\begin{theorem}\label{T:3.1} There is a unique isomorphism of \ $(\mathbf C
L_{-1}\ltimes\tilde{\mathfrak g}_{\leq 0})$-modules
$$
\Xi_W \colon W_{\text{ad}}\otimes V\to U(\tilde{\mathfrak
g})\otimes_{U(\tilde{\mathfrak g}_{\geq 0}+\mathbb C c)} W
$$
such that \  $\Xi_W(w\otimes \mathbf 1)= 1\otimes w$ \ for all
$w\in W$.
\end{theorem}
\begin{proof}
We shall use the identification (\ref{E:3.1}). By definition, a
linear map
$$
\Xi \colon W_{\text{ad}}\otimes U(\tilde{\mathfrak g}_{< 0})\to
U(\tilde{\mathfrak g}_{< 0})\otimes W
$$
is a $\tilde{\mathfrak g}_{<0}$-module homomorphism if
\begin{equation}\label{E:3.7}
\Xi(\,[x,w]\otimes u)+\Xi(w\otimes x_{-1}u)=x_{-1}\Xi(w\otimes u)
\end{equation}
for all $w\in W$, $u\in U(\tilde{\mathfrak g}_{< 0})$ and $x\in
\tilde{\mathfrak g}_{< 0}\mathbf 1\subset V$. Since the universal
enveloping algebra $U(\tilde{\mathfrak g}_{< 0})$ is a quotient of
the tensor algebra $T(\tilde{\mathfrak g}_{< 0})$, we first define
recursively a linear map
$$
\Xi \colon W\otimes T(\tilde{\mathfrak g}_{< 0})\to
U(\tilde{\mathfrak g}_{< 0})\otimes W,
$$
$$
\Xi^s \colon W\otimes T^s(\tilde{\mathfrak g}_{< 0})\to
U_s(\tilde{\mathfrak g}_{< 0})\otimes W,
$$
by setting
$$
\Xi^0(w\otimes \mathbf 1)= 1\otimes w\qquad\text{for all}\qquad
w\in W,
$$
$$
\Xi^{s+1}(w\otimes x_{-1}\otimes u)=x_{-1}\Xi^s(w\otimes
u)-\Xi^s(\,[x,w]\otimes u)
$$
for $s=0,1,2,\dots$. Then we see that $\ker \Xi\supset W\otimes J$
for the ideal $J\subset T(\tilde{\mathfrak g}_{< 0})$ generated by
elements $x_{-1}\otimes y_{-1}-y_{-1}\otimes
x_{-1}-[x_{-1},y_{-1}]$ for $x, y\in \tilde{\mathfrak g}_{<
0}\mathbf 1\subset V$:
\begin{align*}
&\Xi^{s+1}\big(w\otimes (x_{-1}\otimes y_{-1}\otimes u-y_{-1}\otimes x_{-1}\otimes u)\big)\\
&\quad -\Xi^{s}(w\otimes [x_{-1},y_{-1}]\otimes u)\\
&=\big(x_{-1}y_{-1}-y_{-1}x_{-1}-[x_{-1},y_{-1}]\big)\,\Xi^{s-1}(w\otimes u)\\
&\quad +\Xi^{s-1}\big((\,[y,[x,w]]-[x,[y,w]]+[[x,y],w]\,)\otimes  u\big)\\
&=0.
\end{align*}
Here we used the relation $[x_{-1},y_{-1}]\mathbf 1=[x,y]$ for the
inverse of the map (\ref{E:3.6}). By passing to the quotient we
obtain a $\tilde{\mathfrak g}_{<0}$-module homomorphism
\begin{equation}\label{E:3.8}
\Xi_W \colon W_{\text{ad}}\otimes U(\tilde{\mathfrak g}_{< 0})\to
U(\tilde{\mathfrak g}_{< 0})\otimes W,
\end{equation}
\begin{equation}\label{E:3.9}
 \Xi^s_W \colon W_{\text{ad}}\otimes
U_s(\tilde{\mathfrak g}_{< 0})\to U_s(\tilde{\mathfrak g}_{<
0})\otimes W.
\end{equation}
From recursive relations
\begin{equation}\label{E:3.10}
\Xi^{s+1}_W(w\otimes x_{-1} u)=x_{-1}\Xi^s_W(w\otimes
u)-\Xi^s_W(\,[x,w]\otimes u)
\end{equation}
we see by induction that $\Xi_W$ is an isomorphism, uniquely
determined by the requirement $\Xi^0_W(w\otimes \mathbf 1)=
1\otimes w$. By induction (\ref{E:3.10}) implies that the map
$\Xi_W$ intertwines the actions of $L_{-1}$, $L_{0}$ and
${\mathfrak g}(0)$.
\end{proof}

\begin{remarks}
(i) Note that, with the identifications made in (\ref{E:3.8}), by
(\ref{E:3.9}) the map $\Xi_W$ preserves the filtration  inherited
from the filtration $U_s(\tilde{\mathfrak g}_{<0})$, $s\in\mathbb
Z_{\geq 0}$, of the universal enveloping algebra
$U(\tilde{\mathfrak g}_{<0})$.

(ii) Let $\sigma$ be an automorphism of ${\mathfrak g}$. We extend
$\sigma$ to an automorphism of $\tilde{\mathfrak g}$ by
$\sigma(x(n))=(\sigma x)(n)$, $\sigma c=c$, $\sigma d=d$, and, as
well, to $U(\tilde{\mathfrak g})$ and $V$. If $W\subset V$ is
invariant for the action of $\sigma$, than $\sigma$ acts on
$W\otimes V$ and $U(\tilde{\mathfrak
g})\otimes_{U(\tilde{\mathfrak g}_{\geq 0}+\mathbb C c)} W$ in a
natural way. Note that $\Xi_W$ intertwines the action of $\sigma$,
as easily seen from (\ref{E:3.10}).
\end{remarks}

Since  $W\subset V$, we have maps
\begin{equation*}
\begin{array}{ll}
\Psi_W \colon U(\tilde{\mathfrak
g})\otimes_{U(\tilde{\mathfrak g}_{\geq 0}+\mathbb C c)} W
\to V,\qquad &u\otimes w \mapsto uw,\\
{}&{}\\
 \Phi_W \colon W\otimes V\to V,\qquad &u\otimes w \mapsto
u_{-1}w.
\end{array}
\end{equation*}
Note that the map $\Psi_W $ is a homomorphism of $\tilde{\mathfrak
g}$-modules, and that $\Psi_W $ intertwines the actions of
$L_{-1}$, $L_{0}$ and $\sigma$ (if $W$ is invariant for $\sigma$).
Hence $\ker \Psi_W$ is a $\tilde{\mathfrak g}$-module, invariant
for $L_{-1}$, $L_{0}$ and $\sigma$. The following theorem relates
$\ker \Phi_W$ with induced representations of $\tilde{\mathfrak
g}$:

\begin{theorem}\label{T:3.2}
The map $\Phi_W$ is is a homomorphism of $(\mathbf C
L_{-1}\ltimes\tilde{\mathfrak g}_{\leq 0})$-modules and
$\Phi_W=\Psi_W\circ\, \Xi_W$. In particular, $\ker \Phi_W\subset
W_{\text{ad}}\otimes V$ is a $(\mathbf C
L_{-1}\ltimes\tilde{\mathfrak g}_{\leq 0})$-module and
$$
\Xi_W(\ker \Phi_W)=\ker \Psi_W.
$$
\end{theorem}

\begin{remark}
Note that $\Phi_W=\Phi|(W\otimes V)$, and likewise $\Xi_W$ and
$\Psi_W$ are restrictions of $\Xi_V$ and $\Psi_V$ respectively. In
what follows we shall sometimes omit the subscript $W$, and write
only $\Phi$, $\Xi$ and $\Psi$, whenever is clear from a context
which $W\subset V$ is fixed.
\end{remark}

\begin{proof}
First note that for $w\in W$, $u\in V\cong U(\tilde{\mathfrak
g}_{< 0})$ and $x\in \tilde{\mathfrak g}_{< 0}\mathbf 1$ we have
\begin{align*}
&\Phi(w\otimes x_{-1}u)=w_{-1}x_{-1}u=x_{-1}w_{-1}u-[x_{-1},w_{-1}]u\\
&\quad =x_{-1}w_{-1}u-\sum_{n\geq 0}\binom{-1}{\,n}(x_nw)_{-2-n}u\\
&\quad =x_{-1}w_{-1}u-\sum_{n\geq 0}(-1)^n\big(D^{(n+1)}(x_nw)\big)_{-1}u\\
&\quad =x_{-1}w_{-1}u-\big([x, w]\big)_{-1}u=x_{-1}\Phi(w\otimes
u)-\Phi([x, w]\otimes u) .
\end{align*}
Also note that $\Phi(w\otimes \mathbf 1)=w=\Psi\big( \Xi(w\otimes
\mathbf 1)\big)$.  By using the filtration $U_s(\tilde{\mathfrak
g}_{<0})$, $s\in\mathbb Z_{\geq 0}$, we see by induction that
\begin{align*}
&\Psi\big( \Xi(w\otimes x_{-1}u)\big)\\
&\quad =x_{-1}\Psi\big(\Xi(w\otimes
u)\big)-\Psi\big(\Xi([x, w]\otimes u)\big)\\
&\quad =x_{-1}\Phi(w\otimes
u)-\Phi([x, w]\otimes u)\\
&\quad =\Phi(w\otimes x_{-1}u),
\end{align*}
i.e., we see that $\Phi=\Psi\circ \Xi$.
\end{proof}

\section{Sugawara's relations for annihilating fields}

From now on we fix $k\neq -g^\vee$ such that $N(k\Lambda_0)$ is a
reducible $\tilde{\mathfrak g}$-module, and we denote by
$N^1(k\Lambda_0)$ its maximal $\tilde{\mathfrak g}$-submodule. Let
$R\subset N(k\Lambda_0)$ be a nonzero subspace such that:
\begin{itemize}
\item  $R$ is finite-dimensional,
\item  $R$ is invariant for $\tilde{\mathfrak g}_{\geq 0}$ and
$\tilde{\mathfrak g}_{> 0}R=0$,
\item  $R\subset N^1(k\Lambda_0)$.
\end{itemize}
The main example we have in mind is for $k\in\mathbb Z_{> 0}$ and
$R=U(\mathfrak g)x_\theta(-1)^{k+1}\mathbf 1$. Set
$$
\bar R = \mathbb C\text{-span}\{r_n \mid r \in R, n \in \mathbb
Z\},
$$
where $r_n$ denotes a coefficient in the vertex operator $Y(r,z)$.
Then $\bar R$ is a $\tilde{{\mathfrak g}}$-module for the adjoint
action given by the commutator formula
$$
[x_m,r_n]=\sum_{i\geq 0}\binom{m}{i}(x_i r)_{m+n-i},\qquad
x\in{\mathfrak g}, \ r\in R.
$$
From now on we take
$$
W=\bar R\,\mathbf 1\subset V
$$
and we call elements in $\ker \Phi_W$ the relations for
annihilating fields (cf. \cite{P2}). By Theorem~\ref{T:3.2} we may
identify the relations for annihilating fields with elements of
$\ker \Psi_W$, which is easier to study by using the
representation theory of affine Lie algebras. The first step is to
introduce Sugawara's relations for annihilating fields as elements
in
$$
N= U(\tilde{\mathfrak g})\otimes_{U(\tilde{\mathfrak g}_{\geq
0}+\mathbb C c)} W\cong U(\tilde{\mathfrak g}_{< 0})\otimes W.
$$

Let $\{x^i\}_{i\in I}$ and $\{y^i\}_{i\in I}$ be dual bases in
$\mathfrak g$. For $r\in R$ we define Sugawara's relation
\begin{equation}\label{E:4.1}
q_r=\frac{1}{k+g^\vee}\sum_{i\in I}x^i(-1)\otimes y^i(0)r-1\otimes
Dr
\end{equation}
as an element of $U(\tilde{\mathfrak
g})\otimes_{U(\tilde{\mathfrak g}_{\geq 0}+\mathbb C c)} W$. As in
the case of Casimir operator, Sugawara's relation $q_r$ does not
depend on a choice of dual bases $\{x^i\}_{i\in I}$ and
$\{y^i\}_{i\in I}$. Since on $V$ we have $D=L_{-1}$, Sugawara's
construction gives (cf. \cite{MP1}, Chapter 3)
\begin{align*}
L_{-1}r&=\frac{1}{2(k+g^\vee)}\sum_{i\in I}\big(x^i(-1)
y^i(0)+y^i(-1) x^i(0)\big)r\\
&=\frac{1}{k+g^\vee}\sum_{i\in I}x^i(-1) y^i(0)r,
\end{align*}
and we have the following proposition (cf. \cite{MP1}, Chapter 8):
\begin{proposition}\label{P:4.1}
(i) $q_r$ is an element of \ $\ker \Psi_W$.

(ii) $r\mapsto q_r$ is a $\mathfrak g$-module homomorphism from
$R$ into $\ker \Psi_W$.

(iii) $x(i)q_r=0$ for all $x\in\mathfrak g$ and $i>0$.
\end{proposition}
Let us denote the set of all Sugawara's relations (\ref{E:4.1}) by
$$
Q_{\text{Sugawara}}=\{q_r\mid r\in R\}\subset\ker \Psi_W.
$$
Our assumptions imply that $L_1R=0$, $L_0R\subset R$
 and $R\cap DR=0$. So $\mathfrak g$-module $Q_{\text{Sugawara}}$ is
isomorphic to $R$, and by Proposition~\ref{P:4.1}(iii) for each
maximal vector $r$ in $R$ the corresponding $q_r$ is a singular
vector in the $\tilde{\mathfrak g}$-module $\ker \Psi_W$.

Our assumptions imply that \ $\bar R\,\mathbf 1=\coprod_{i\geq 0}
D^iR $ \ and we set
\begin{equation*}
W_n=\coprod_{i=0}^n D^iR,\qquad N_n=U(\tilde{\mathfrak
g})\otimes_{U(\tilde{\mathfrak g}_{\geq 0}+\mathbb C c)} W_n.
\end{equation*}
So we have a filtration of  $N$  by $\tilde{\mathfrak g}$-modules
$N_n$, $n\in\mathbb Z_{\geq 0}$, and $N_0$ is a (sum of)
generalized Verma module(s). We denote the restrictions
$\Psi_W|N_n$ simply by $\Psi_n$. In particular, we have the
$\tilde{\mathfrak g}$-module homomorphism
\begin{equation}\label{E:4.2}
\Psi_0  \colon U(\tilde{\mathfrak g})\otimes_{U(\tilde{\mathfrak
g}_{\geq 0}+\mathbb C c)} R \to N^1(k\Lambda_0),\qquad u\otimes w
\mapsto uw.
\end{equation}
\begin{proposition}\label{P:4.2} As a $(\mathbf C
L_{-1}\ltimes\tilde{\mathfrak g}_{\leq 0})$-module  $\ker \Psi_W$
is generated by
$$
\ker\Psi_0+Q_{\text{Sugawara}}.
$$
\end{proposition}
\begin{proof}
If $v\in N_0\cap\ker \Psi_W$, then $v\in\ker\Psi_0$. Now assume
that $v\in N_n\cap\ker \Psi_W$ for $n>0$. Then $v$ is a sum of
elements of the form
$$
v=\sum_{u,r}u\otimes D^nr+w,\qquad u\in U(\tilde{\mathfrak g}_{<
0}),\ r\in R,\ w\in N_{n-1},
$$
and to each summand $u\otimes D^nr$ add
$$
uD^{n-1}q_r=uD^{n-1}\left(\frac{1}{k+g^\vee}\sum_{i\in
I}x^i(-1)\otimes y^i(0)r\right)-u\otimes  D^nr\in\ker \Psi_W.
$$
Then $v+\sum uD^{n-1}q_r\in N_{n-1}\cap\ker \Psi_W$, and in a
finite number of steps we get
$$
v\in\ker\Psi_0+\sum_{i=0}^{n-1}U(\tilde{\mathfrak g}_{< 0})D^i
Q_{\text{Sugawara}}.
$$
\end{proof}
\begin{remark}
By using (\ref{E:3.4}) and (\ref{E:3.7}) it is easy to see that
for $r\in R$ the corresponding Sugawara's relation $\Xi^{-1}(q_r)$
in $\ker\Phi_W\subset W\otimes V$ is
\begin{equation}\label{E:4.3}
\Xi^{-1}(q_r)=\frac{1}{k+g^\vee}\sum_{i\in I}y^i(0)r \otimes
x^i(-1)\mathbf 1 + \frac{1}{k+g^\vee}\,D\Omega r\otimes\mathbf
1-Dr\otimes\mathbf 1,
\end{equation}
where $\Omega=\sum_{i\in I}x^i(0)y^i(0)$ is the
Casimir operator for $\mathfrak g(0)$.
\end{remark}

 \section{Generators of $\ker\Psi_0$ in the case of standard modules}

From now on we consider the case of standard $\tilde{\mathfrak
g}$-modules $L(k\Lambda_0)$: we take $k\in\mathbb Z_{> 0}$ and
$R=U(\mathfrak g)x_\theta(-1)^{k+1}\mathbf 1$, where
$x_\theta\in\mathfrak g$ is a root vector for the maximal root
$\theta$. In this case the maximal $\tilde{\mathfrak g}$-submodule
$N^1(k\Lambda_0)$ of $N(k\Lambda_0)$ is generated by the singular
vector $x_\theta(-1)^{k+1}\mathbf 1$, so that $N_0$ is a
generalized Verma module and the map $\Psi_0$ defined by
(\ref{E:4.2}) is surjective (cf. \cite{MP1}). Hence we have the
exact sequence of $\tilde{\mathfrak g}$-modules
$$
N_0 \xrightarrow{\Psi_0} N(k\Lambda_0)\rightarrow
L(k\Lambda_0)\rightarrow 0.
$$
We shall use Garland-Lepowsky's resolution \cite{GL} of a standard
module, in terms of generalized Verma modules, to determine
generators of $\ker\Psi_0$.

Denote by $W(\tilde{\mathfrak g})$ the Weyl group generated by
simple reflections $r_0, r_1,\dots, r_\ell$. Note that for
$i\in\{1,\dots,\ell\}$
$$
r_0r_i\neq r_ir_0 \quad\text{if and only if}\quad
\langle\alpha_0,\alpha_i^\vee\rangle\neq 0,
$$
and that for $A_\ell^{(1)}$, $\ell\geq 2$, there are exactly two
such $i$, for all the other untwisted affine Lie algebras
$\tilde{\mathfrak g}$ there is exactly one such $i$, corresponding
to $\alpha_i$ connected with $\alpha_0$ in a Dynkin diagram. With
the usual notation, for $w\in W(\tilde{\mathfrak g})$ and a weight
$\lambda$ write
$$
w\cdot \lambda=w(\lambda+\rho)-\rho.
$$
\begin{lemma}\label{L:5.1} For $\tilde{\mathfrak g}$ of the type
$A_\ell^{(1)}$, $\ell\geq 2$, the $\tilde{\mathfrak g}$-module
$\ker\Psi_0$ is generated by exactly two singular vectors of
weights $r_0r_i\cdot k\Lambda_0$,
$\langle\alpha_0,\alpha_i^\vee\rangle\neq 0$. For all the other
untwisted affine Lie algebras $\tilde{\mathfrak g}$ the
$\tilde{\mathfrak g}$-module $\ker\Psi_0$ is generated by exactly
one singular vector of weight $r_0r_i\cdot k\Lambda_0$,
$\langle\alpha_0,\alpha_i^\vee\rangle\neq 0$.
\end{lemma}
\begin{proof}
A generalized Verma $\tilde{\mathfrak g}$-module $N(\Lambda)$,
with $\Lambda(\alpha_i^\vee)\in\mathbb Z_{\geq 0}$ for
$i\in\{1,\dots,\ell\}$ and $\Lambda(c)=k$, may be defined (cf.
\cite{Le}) as an induced $\tilde{\mathfrak g}$-module
$U(\tilde{\mathfrak g})\otimes_{U(\tilde{\mathfrak g}_{\geq
0}+\mathbb C c)}L $, where $L=L(\Lambda|\mathfrak h)$ is the
irreducible finite-dimensional ${\mathfrak g}$-module with the
highest weight $\Lambda|\mathfrak h$, on which $\tilde{\mathfrak
g}_{> 0}$ acts trivially and $c$ as the multiplication by $k$. In
terms of Verma modules the generalized Verma module $N(\Lambda)$
may be written as a quotient
\begin{equation}\label{E:5.1}
N(\Lambda)=M(\Lambda)\big/\big(\sum_{i=1}^\ell
M(r_i\cdot\Lambda)\big).
\end{equation}
By Theorem~8.7 in \cite{GL} (cf. Theorem~5.1 in \cite{Le}) we have
Garland-Lepowsky's resolution of $L(k\Lambda_0)$ in terms of
generalized Verma modules
$$
\dots \rightarrow E_2 \rightarrow E_1 \rightarrow E_0\rightarrow
L(k\Lambda_0)\rightarrow 0.
$$
In our setting one gets $E_0=N(k\Lambda_0)$ and $E_1=N(r_0\cdot
k\Lambda_0)=N_0$. In the case when $\tilde{\mathfrak g}$ is of the
type $A_\ell^{(1)}$, $\ell\geq 2$, the $\tilde{\mathfrak
g}$-module $E_2$ has a filtration $0=V_0\subset V_1\subset
V_2=E_2$ such that $V_1/V_0$ and $V_2/V_1$ are exactly two
generalized Verma modules $N(r_0r_i\cdot k\Lambda_0)$,
$\langle\alpha_0,\alpha_i^\vee\rangle\neq 0$. Hence $\ker\Psi_0$
is generated by one singular vector and one subsingular vector,
their two weights being $r_0r_i\cdot k\Lambda_0$,
$\langle\alpha_0,\alpha_i^\vee\rangle\neq 0$. Since $R$ and
$\ker\Psi_0$ are invariant for the action of the Dynkin diagram
automorphism $\sigma$, which interchanges these two weights, both
generating vectors must be singular vectors. For all the other
untwisted affine Lie algebras $\tilde{\mathfrak g}$ one gets that
$E_2=N(r_0r_i\cdot k\Lambda_0)$, and hence the $\tilde{\mathfrak
g}$-module $\ker\Psi_0$ is generated by exactly one singular
vector of weight $r_0r_i\cdot k\Lambda_0$,
$\langle\alpha_0,\alpha_i^\vee\rangle\neq 0$.
\end{proof}
\begin{remark}
One can also prove Lemma~\ref{L:5.1} by using the BGG type
resolution of a standard module in terms of Verma modules, due to
A.~Rocha-Caridi and N.~R.~Wallach \cite{RW}.
\end{remark}

\begin{remark}
Singular vectors in Verma module $M(r_0\cdot k\Lambda_0)$ can be
computed by the Malikov-Feigin-Fuchs formula \cite{MFF}, and by
passing to the quotient (\ref{E:5.1}) we may find singular vectors
in $N_0=N(r_0\cdot k\Lambda_0)$. In the case $A_1^{(1)}$ level
$k=1$ the singular vector of weight $\lambda=r_0r_1\cdot
\Lambda_0$ in $M(r_0\cdot \Lambda_0)$ is given by the formula
$$
f_0^{4}f_1^{1}f_0^{-2}v_{r_0\cdot \Lambda_0},
$$
with the usual notation for generators of Kac-Moody Lie algebras.
This formula gives the singular vector
$$
\big(x_{\theta}(-1)^2x_{-\theta}(0)+
2x_{\theta}(-1)\theta^\vee(-1)-6x_{\theta}(-2)\big)v_{r_0\cdot
\Lambda_0}
$$
in $M(r_0\cdot \Lambda_0)$. So we have a generator of
$\ker\Psi_0\subset N_0=N(r_0\cdot \Lambda_0)$
\begin{equation}\label{E:5.2}
\big(2x_{\theta}(-1)\theta^\vee(-1)-6x_{\theta}(-2)\big)\otimes
x_{\theta}(-1)^2\mathbf 1+ x_{\theta}(-1)^2\otimes
x_{-\theta}(0)x_{\theta}(-1)^2\mathbf 1.
\end{equation}
We shall determine singular vectors in $\ker\Psi_0$ in another
way.
\end{remark}

Let us denote by $\alpha_*$ and $r_*$ the elements $\alpha_i$ and
$r_i$ for which $\langle\alpha_0,\alpha_i^\vee\rangle\neq 0$. Then
\begin{equation}\label{E:5.3}
r_0r_*\cdot k\Lambda_0=k\Lambda_0-\alpha_*
-\big(k+1-\langle\alpha_*,\alpha_0^\vee\rangle\big)\alpha_0.
\end{equation}
Note that $\langle\alpha_*,\alpha_0^\vee\rangle=-2$ for
$A_1^{(1)}$, and $\langle\alpha_*,\alpha_0^\vee\rangle=-1$ for all
the other untwisted affine Lie algebras $\tilde{\mathfrak g}$. For
this reason the case $A_1^{(1)}$ is somewhat different, so for now
assume that $\tilde{\mathfrak g}\not\cong{\mathfrak sl}(2,\mathbb
C)\,\widetilde{}$\,. Then

\begin{lemma}\label{L:5.2}
$\theta-\alpha_*$ is a root of ${\mathfrak g}$ and \
$2\theta-\alpha_*$ is not a root of ${\mathfrak g}$.
\end{lemma}
\begin{proof}
For the type $C_\ell$, with the usual notation,
$\theta=2\varepsilon_1$ and
$\alpha_*=\alpha_1=\varepsilon_1-\varepsilon_2$, and the statement
is clear. For all the other types
$\langle\theta,\alpha_*^\vee\rangle=1$, so
$r_*\,\theta=\theta-\alpha_*$ is a positive root of ${\mathfrak
g}$, and the statement is clear.
\end{proof}
Then $x_{\theta-\alpha_*}=[x_{-\alpha_*},x_{\theta}]$ is a root
vector in ${\mathfrak g}_{\theta-\alpha_*}$, and again by
Lemma~\ref{L:5.2},
$$
[x_{\theta-\alpha_*},x_{\theta}]=0.
$$
Since
$$
x_{-\alpha_*}(0)x_\theta(-1)^{k+1}\mathbf
1=(k+1)x_{\theta-\alpha_*}(-1)x_\theta(-1)^{k}\mathbf 1\in R,
$$
we obviously have
$$
x_{\theta-\alpha_*}(-1)\otimes x_\theta(-1)^{k+1}\mathbf 1\mathbf
-x_{\theta}(-1)\otimes
x_{\theta-\alpha_*}(-1)x_\theta(-1)^{k}\mathbf 1\in\ker\Psi_0.
$$
By (\ref{E:5.3}) this vector is of weight $r_0r_*\cdot
k\Lambda_0$, and hence by Lemma~\ref{L:5.1} it must be a singular
vector. So we have proved
\begin{proposition}\label{P:5.3}
Let $\tilde{\mathfrak g}\not\cong{\mathfrak sl}(2,\mathbb
C)\,\widetilde{}$\,  be an untwisted affine Lie algebra. Then
$\ker\Psi_0$ is generated by the singular vector(s)
$$
x_{\theta-\alpha_*}(-1)\otimes x_\theta(-1)^{k+1}\mathbf 1\mathbf
-x_{\theta}(-1)\otimes
x_{\theta-\alpha_*}(-1)x_\theta(-1)^{k}\mathbf 1,\quad
\langle\alpha_0,\alpha_*^\vee\rangle\neq 0.
$$
\end{proposition}
\bigskip

Now let $\tilde{\mathfrak g}={\mathfrak sl}(2,\mathbb
C)\,\widetilde{}$\,. We have the Sugawara singular vector
\begin{align*}
q_{(k+1)\theta}=&\tfrac{1}{k+2}\big(x_{\theta}(-1)\otimes
x_{-\theta}(0)x_{\theta}(-1)^{k+1}\mathbf 1+\tfrac12
\theta^\vee(-1)\otimes
\theta^\vee(0)x_{\theta}(-1)^{k+1}\mathbf 1\big)\\
&-1\otimes Dx_{\theta}(-1)^{k+1}\mathbf 1,
\end{align*}
and the obvious relation
$$
q_{(k+2)\theta}=x_\theta(-2)\otimes x_\theta(-1)^{k+1}\mathbf 1-
\tfrac{1}{k+1}\,x_\theta(-1)\otimes Dx_\theta(-1)^{k+1}\mathbf 1
$$
(cf. \cite{MP1}). Then obviously
$$
(k+1)q_{(k+2)\theta}-x_\theta(-1)q_{(k+1)\theta}\in\ker\Psi_0
$$
is a nonzero vector of weight $r_0r_1\cdot
k\Lambda_0=k\Lambda_0-\alpha_1 -(k+3)\alpha_0$, and hence by
Lemma~\ref{L:5.1} it must be a singular vector, proportional to
(\ref{E:5.2}) in the case $k=1$.

By combining Theorem~\ref{T:3.2}, Propositions~\ref{P:4.2},
\ref{P:5.3}, and (\ref{E:4.3}), we get:
\begin{theorem}\label{T:5.4}
Let $\tilde{\mathfrak g}\not\cong{\mathfrak sl}(2,\mathbb
C)\,\widetilde{}$\,  be an untwisted affine Lie algebra. Then the
$(\mathbf C L_{-1}\ltimes\tilde{\mathfrak g}_{\leq 0})$-module
$\ker \Phi_W$ is generated by vectors
$$
 x_\theta(-1)^{k+1}\mathbf 1 \otimes x_{\theta-\alpha_*}(-1)\mathbf 1
- x_{\theta-\alpha_*}(-1)x_\theta(-1)^{k}\mathbf 1\otimes
x_{\theta}(-1)\mathbf 1,\quad
\langle\alpha_0,\alpha_*^\vee\rangle\neq 0,
$$
$$
\tfrac{1}{k+g^\vee}\sum_{i\in I}y^i(0)x_\theta(-1)^{k+1}\mathbf 1
\otimes x^i(-1)\mathbf 1 + L_{-1}\big(\tfrac{1}{k+g^\vee}\,\Omega
-1\big)x_\theta(-1)^{k+1}\mathbf 1 \otimes\mathbf 1.
$$
\end{theorem}

\section{Relation \ $x_\theta(z)\tfrac{d}{dz}(x_\theta(z)^{k+1})=
(k+1)x_\theta(z)^{k+1}\tfrac{d}{dz}x_\theta(z)$}
\bigskip

We have the obvious relation in $\ker\Psi_W$
\begin{equation}\label{E:6.1}
q_{(k+2)\theta}=x_{\theta}(-2)\otimes x_\theta(-1)^{k+1}\mathbf
1\mathbf -x_{\theta}(-1)\otimes
x_{\theta}(-2)x_\theta(-1)^{k}\mathbf 1,
\end{equation}
and the corresponding relation in $\ker\Phi_W$ is
$$
 \Xi^{-1}(q_{(k+2)\theta})=x_\theta(-1)^{k+1}\mathbf 1 \otimes x_{\theta}(-2)\mathbf 1
- x_{\theta}(-2)x_\theta(-1)^{k}\mathbf 1\otimes
x_{\theta}(-1)\mathbf 1.
$$
The vertex operator associated to the vector
$\Xi^{-1}(q_{(k+2)\theta})$ in the vertex operator algebra
$V\otimes V$ is
\begin{equation}\label{E:6.2}
x_\theta(z)^{k+1}\otimes\tfrac{d}{dz}x_\theta(z)-
\tfrac{1}{k+1}\,\tfrac{d}{dz}(x_\theta(z)^{k+1})\otimes
x_\theta(z),
\end{equation}
and by taking normal order products we get for the annihilating
field $x_\theta(z)^{k+1}$ the obvious relation
\begin{equation}\label{E:6.3}
x_\theta(z)^{k+1}\tfrac{d}{dz}x_\theta(z)-
\tfrac{1}{k+1}\,\tfrac{d}{dz}(x_\theta(z)^{k+1}) x_\theta(z)=0.
\end{equation}

The relation (\ref{E:6.3}) plays a key role in constructions
\cite{MP1}, \cite{MP2}, \cite{S} of combinatorial bases of
standard ${\mathfrak sl}(2,\mathbb C)\,\widetilde{}$\, and
${\mathfrak sl}(3,\mathbb C)\,\widetilde{}$\,-modules, and these
examples showed that all necessary relations for annihilating
fields of standard ${\mathfrak sl}(2,\mathbb C)\,\widetilde{}$\,
and ${\mathfrak sl}(3,\mathbb C)\,\widetilde{}$\,-modules can be
derived from (\ref{E:6.3}), or (\ref{E:6.2}), to be a bit more
precise. In this section we show that for a general untwisted
affine Lie algebra $\tilde{\mathfrak g}\not\cong{\mathfrak
sl}(2,\mathbb C)\,\widetilde{}$ \  we can obtain from
(\ref{E:6.1}), by using the action of $\tilde{\mathfrak g}_{\geq
0}$ on $\ker\Psi_W$, both the singular vector(s)
$$
q_{(k+2)\theta-\alpha_*}=x_{\theta-\alpha_*}(-1)\otimes
x_\theta(-1)^{k+1}\mathbf 1\mathbf -x_{\theta}(-1)\otimes
x_{\theta-\alpha_*}(-1)x_\theta(-1)^{k}\mathbf 1
$$
in $\ker\Psi_0$ and the Sugawara singular vector
$$
q_{(k+1)\theta}=\tfrac{1}{k+g^\vee}\sum_{i\in I}x^i(-1)\otimes
y^i(0)x_\theta(-1)^{k+1}\mathbf 1-1\otimes
Dx_\theta(-1)^{k+1}\mathbf 1.
$$
\begin{lemma}\label{L:6.1} Let $\Omega$ be the Casimir
operator for ${\mathfrak g}\not\cong{\mathfrak sl}(2,\mathbb C)$
and $\lambda=(k+2)\theta-\alpha_*$. Then
\begin{align*}
&q_{(k+2)\theta-\alpha_*}=x_{-\alpha_*}(1)q_{(k+2)\theta},\\
&q_{(k+1)\theta}=\tfrac{k+1}{2(k+2)(k+g^\vee)}\big(\Omega-(\lambda+2\rho,\lambda)\big)
x_{-\theta}(1)q_{(k+2)\theta}.
\end{align*}
\end{lemma}
\begin{proof}
Since $x_\theta(-1)^{k+1}\mathbf 1$ is a singular vector, and
$[x_{\theta-\alpha_*},x_{\theta}]=0$, the first equality is
obvious. To prove the second equality, first note that
\begin{align*}
(k+1)x_{-\theta}(1)q_{(k+2)\theta}=&(k+2)\otimes D
x_\theta(-1)^{k+1}\mathbf 1\\
-&x_\theta(-1)\otimes
x_{-\theta}x_\theta(-1)^{k+1}\mathbf 1\\
-&(k+1)\theta^\vee(-1)\otimes
 x_\theta(-1)^{k+1}\mathbf 1,
\end{align*}
so that
$$
q=(k+1)x_{-\theta}(1)q_{(k+2)\theta}+(k+2)q_{(k+1)\theta}
$$
is in $\ker\Psi_0$. Since $q$ is of degree $-k-2$, the same as the
degree of singular vector(s) in $\ker\Psi_0$, $q$ must be an
element of a $\mathfrak g$-module generated by vector(s)
$q_{(k+2)\theta-\alpha_*}$. Hence $\Omega
q=(\lambda+2\rho,\lambda) q$  and
\begin{align*}
&(k+1)\big(\Omega -(\lambda+2\rho,\lambda)\big)x_{-\theta}(1)q_{(k+2)\theta}\\
&=(k+2)\big((\lambda+2\rho,\lambda)-\Omega\big)q_{(k+1)\theta}\\
&=(k+2)\cdot \big((\lambda+2\rho,\lambda)-
((k+1)\theta+2\rho,(k+1)\theta)\big)\cdot q_{(k+1)\theta}\\
&=(k+2)\cdot 2(k+g^\vee)\cdot q_{(k+1)\theta}
\end{align*}
gives the second formula.
\end{proof}
It is easy to see that for ${\mathfrak g}={\mathfrak sl}(2,\mathbb
C)\,\widetilde{}$ \ we have
$$
(k+1)x_{-\theta}(1)q_{(k+2)\theta}+(k+2)q_{(k+1)\theta}=0.
$$
Hence, by combining Proposition~\ref{P:4.2},
Proposition~\ref{P:5.3} and Lemma~\ref{L:6.1} we get:
\begin{theorem}\label{T:6.2}
Let $\tilde{\mathfrak g}$  be an untwisted affine Lie algebra.
Then the $(\mathbf C L_{-1}\ltimes\tilde{\mathfrak g})$-module
$\ker \Psi_W$ is generated by the vector $q_{(k+2)\theta}$.
\end{theorem}
\begin{remarks}
(i) Although Theorem~\ref{T:5.4} describes generators of relations
for annihilating fields, for combinatorial applications (cf.
\cite{P2}, Theorem~2.12) Theorem~\ref{T:6.2}, in conjunction with
Theorem~\ref{T:3.2}, seems to be better suited. Namely, what one
needs is a description of
$$
(\ker\Phi_W)_s=(V\otimes V)_s\cap\ker\Phi_W,
$$
where $(V\otimes V)_s$, $s\in\mathbb Z_{\geq 0}$, is the natural
filtration of $V\otimes V$ inherited from the filtration
$U_s(\tilde{\mathfrak g}_{<0})$, $s\in\mathbb Z_{\geq 0}$. For the
corresponding filtration $(\ker\Psi_W)_s$, $s\in\mathbb Z_{\geq
0}$, of $\ker\Psi_W$ the subspace $(\ker\Psi_W)_s$ is not
preserved by the action of $\tilde{\mathfrak g}_{< 0}$, but it is
preserved by the action of $\tilde{\mathfrak g}_{\geq 0}$.
Theorem~\ref{T:6.2} makes one think that $D^n q_{(k+2)\theta}$,
$n\in\mathbb Z_{\geq 0}$, might be generators of
$(\ker\Psi_W)_{k+2}$ for the action of $\tilde{\mathfrak g}_{\geq
0}$\,?

(ii) In the case of level $k=1$ standard ${\mathfrak sl}(3,\mathbb
C)\,\widetilde{}$\,-modules  one needs only the relations
$q_{3\theta}$, $q_{3\theta-\alpha_1}$, $q_{3\theta-\alpha_2}$,
$q_{2\theta}\in (\ker\Psi_W)_3$ for a construction of
combinatorial bases \cite{MP2}, \cite{S}.
\end{remarks}

\end{document}